\documentclass{ro}
\usepackage{epstopdf}
\usepackage{subfigure}
\usepackage{url}
\usepackage{mathtools}
\usepackage{booktabs}
\usepackage{color}
\usepackage{amsthm}
\theoremstyle{remark}

\begin{document}
\title{Simultaneous estimation of multiple discrete unimodal distributions under stochastic order constraints}
%
\author{Y. Yoshida}
\address{School of Engineering, Stanford University, Stanford, CA 94305, USA}

\author{N. Sukegawa}
\address{Faculty of Science and Engineering, Hosei University, Tokyo, Japan}

\author{J. Iwanaga}
\address{Erdos Inc., Kanagawa, Japan; The University of Electro-Communications, Tokyo, Japan}
\date{31 July 2026}
\begin{abstract}
We study the problem of estimating multiple discrete unimodal distributions, motivated by search behavior analysis on a real-world platform. To incorporate prior knowledge of precedence relations among distributions, we impose stochastic order constraints and formulate the estimation task as a mixed-integer convex quadratic optimization problem. Experiments on both synthetic and real datasets show that the proposed method reduces the Jensen–Shannon divergence by 2.2\% on average (up to 6.3\%) when the sample size is small, while performing comparably to existing methods when sufficient data are available.
\end{abstract}
%
%
\maketitle
\section{Introduction}
\label{sec:intro}
In this paper, we present a mixed-integer convex optimization model to simultaneously estimate multiple discrete unimodal distributions under stochastic order constraints, motivated by a real-world search behavior analysis. 

\subsection{Background}
Maternal mental health has become a serious concern, especially during the Covid-19 pandemic~\cite{davenport2020moms}. 
The WHO reports that about 10\% of pregnant women and 13\% of postpartum women experience mental disorders, mainly depression~\cite{who}. In this context, Connehito Inc. operates Mamari, an information platform for pregnancy, childbirth, and childcare, featuring user-to-user Q\&A. By 2024, Mamari is projected to have 3.5 million users, with 4 million monthly searches and 1.1 million monthly Q\&A posts.

\begin{figure}[tb]
\includegraphics[clip, width=\textwidth]{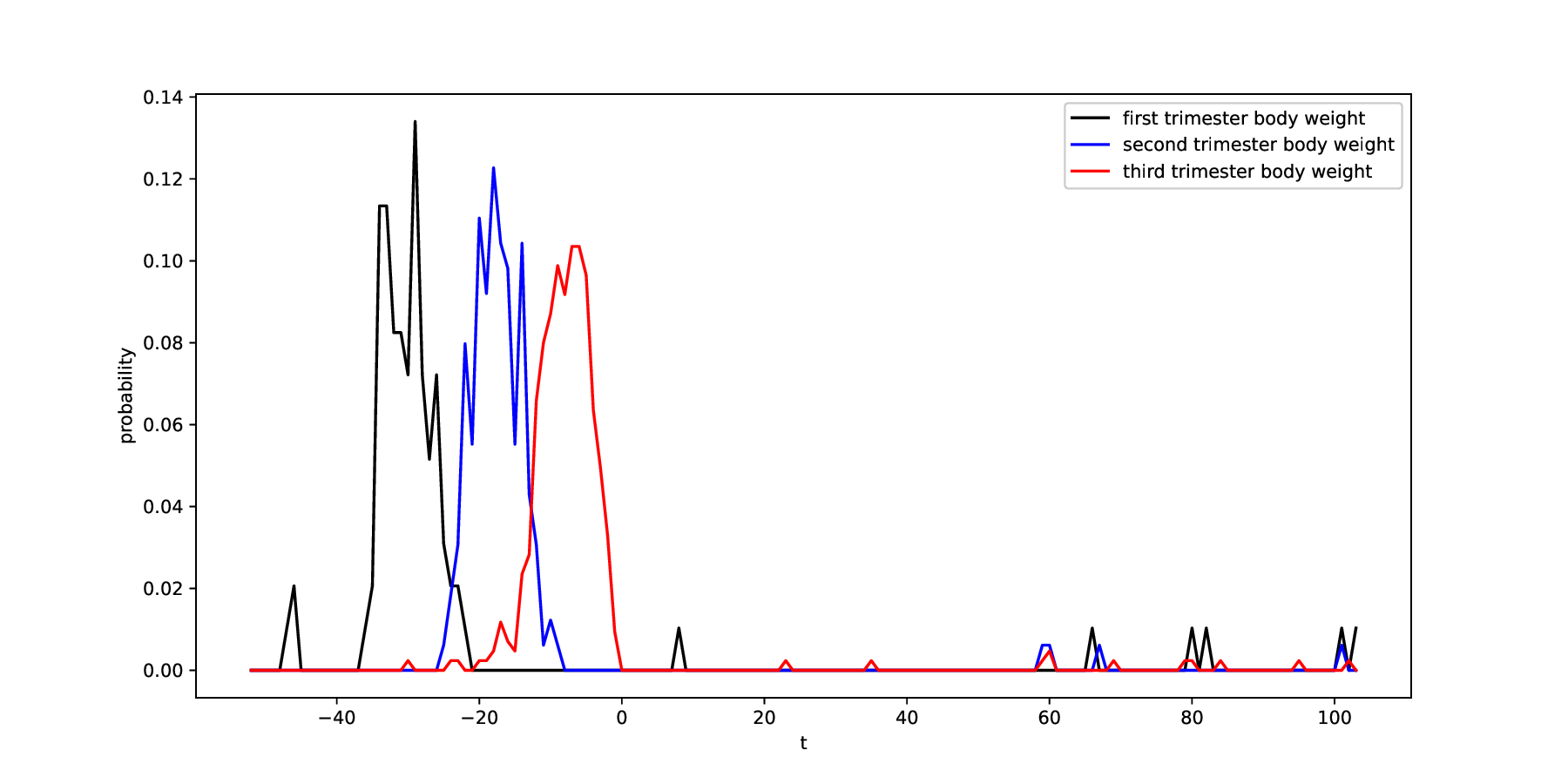}
\caption{Three search timing distributions for keywords containing 
\texttt{body weight} on Mamari}
\label{fig:bw}
\end{figure}

Understanding user interest through search queries is a fundamental challenge in information retrieval~\cite{qiu2006automatic}. 
As reported in \cite{iwanaga2022analysis}, search timing for typical Mamari keywords depends on children’s ages.  
Specifically, the \emph{search timing distributions} tended to be unimodal, as depicted in Figure~\ref{fig:bw}. 
In this figure, the black line, for example, shows the normalized histogram of the number of users who searched for \texttt{first\_trimester\_body\_weight} per week, with the $x$-axis indicating the age of the children. 
This figure suggests that most users in the first trimester are interested in their body weight approximately $30$ weeks before the expected date of birth, which is indicated by zero on the $x$-axis. 

\subsection{Motivation}
In \cite{iwanaga2022analysis}, unimodal regression~\cite{stout2000optimal} outperformed baseline methods in estimating search timing distributions for typical keywords. The estimation error, however, increases for multi-term keywords, such as \texttt{first\_trimester\_body\_weight}, owing to smaller sample sizes.

To mitigate this limitation, we exploit prior knowledge of search timing distributions. By definition, the distribution of \texttt{first trimester\_body weight} should \emph{precede} that of \texttt{second trimester\_body weight}. Figure~\ref{fig:bw} validates this assumption and suggests that samples from \texttt{first trimester} may also improve estimation for \texttt{second trimester}. 

\subsection{Our contribution}
In this study, we propose a model that incorporates prior knowledge of \emph{precedence relations} among search timing distributions to improve estimation accuracy under limited samples. The main contributions are:
\begin{enumerate}
\item We formalize precedence relations via stochastic order~\cite{shaked2007stochastic} and show that the resulting estimation problem reduces to a mixed-integer convex quadratic program solvable by standard solvers such as Gurobi.
\item Using real search history data from Mamari, we demonstrate that the proposed model substantially reduces estimation error compared with baseline methods.
\end{enumerate}

The remainder of this paper is organized as follows. Section~\ref{sec:rw} reviews related work. Section~\ref{sec:mm} presents the proposed model and its formulation as a mixed-integer convex quadratic program. Section~\ref{sec:poc} provides a proof of concept, followed by numerical results in Section~\ref{sec:ne}. Section~\ref{sec:conc} concludes the paper.

\section{Related work}
\label{sec:rw}

Distribution estimation has been widely studied~\cite{silverman1986density,scott2015multivariate}, with most work focusing on single distributions or the independent estimation of multiple distributions. 
In contrast, this work considers the simultaneous estimation of multiple discrete distributions subject to cross-distribution structural constraints.

A common approach to distribution estimation is nonparametric modeling, such as kernel density estimation~\cite{chen2017tutorial,fadda1998density,sheather2004density}. While highly flexible, nonparametric models are prone to overfitting, motivating the incorporation of application-specific prior knowledge~\cite{von2021informed,wu2022survey,diligenti2017integrating}. Typical prior knowledge includes shape constraints such as unimodality~\cite{turnbull2014unimodal,hall2002unimodal}, log-concavity~\cite{samworth2018recent}, and $k$-monotonicity~\cite{balabdaoui2007estimation}. Monotonicity leads to isotonic regression~\cite{barlow1972isotonic}, for which efficient algorithms exist~\cite{kalai2009isotron,tibshirani2011nearly}. Unimodal estimation has also been studied in discrete settings, notably through the optimal pyramid~\cite{chun2003efficient,chun2006efficient}, which focuses on single distributions. 
However, existing frameworks do not handle stochastic order constraints across multiple distributions or jointly enforce inter-distribution ordering and shape constraints within a unified optimization model.

Stochastic orders play an important role in probability theory and statistics as tools for comparing stochastic models and deriving probabilistic inequalities~\cite{bergmann1991stochastic,mosler2012stochastic}. 
Most existing work on stochastic dominance focuses on testing and inference rather than estimation~\cite{lando2025new,barrett2003consistent}. 
There are also studies on optimization under stochastic order constraints~\cite{dentcheva2016augmented,dentcheva2016two,dentcheva2009optimization,muller2006stochastic}, including analytical results for maximum likelihood estimation of two distributions under a stochastic order constraint~\cite{el2005inferences}. 
However, these approaches are typically restricted to pairwise ordering and specific estimation settings. In contrast, we propose a mixed-integer optimization framework for the simultaneous estimation of multiple discrete distributions that jointly enforces unimodality, stochastic order constraints across more than two distributions, and additional structural constraints such as bounded support, while directly minimizing deviation from empirical distributions.

\section{Our model}
\label{sec:mm}
We review stochastic order for continuous and discrete distributions~\cite{bauerle2006stochastic,gupta2010convex,liu2003testing} and formulate the estimation problem as a mixed-integer convex quadratic optimization problem. 

\subsection{Stochastic order}
For two random variables $X_1$ and $X_2$, we say that $X_1$ is \emph{stochastically smaller than} $X_2$ and write $X_1 \le_{\rm st} X_2$ if $E[f(X_1)] \le E[f(X_2)]$ for all increasing functions $f$. 
It is known that $X_1 \le_{\rm st} X_2$ holds if and only if
$F_{X_1}(t) \ge F_{X_2}(t)$ holds for any $t \in \mathbb{R}$, where $F_X$ denotes the distribution function of $X$. 

From the above characterization, for any two discrete distributions $P_1$ and $P_2$ with a common support $T=\{l,\ldots,u\}$, we have $P_1 \le_{\rm st} P_2$ if 
\begin{align}
\label{st_d}
\sum_{i \le t} p_{1i} \ge \sum_{i \le t} p_{2i}
\end{align}
holds for any $t \in T$, where $p_{1i}$ and $p_{2i}$ denote the $i$th components of $P_1$ and $P_2$, respectively. 

\subsection{Unimodal regression}
\label{sec:unimodal}
A discrete distribution $P$ is \emph{unimodal} if $p_i \le p_{i+1}$ for $i < t$ and $p_i \ge p_{i+1}$ for $t \le i$ for some $t\in T$.
The unimodal regression fits a unimodal distribution $X$ to an empirical distribution $P$, which can be formulated as the following 0-1 mixed integer convex programming problem:
\begin{align*}
\begin{array}{lllll}
R(P):~&\mbox{min.}    &d(X,P)\\\\
&\mbox{s. t.}   &\displaystyle \sum_{i\in T} x_i = 1 &&\mbox{(C1)}\\
&&x_i \in [0,1]     &(i\in T) &\mbox{(C2)}\\
&&x_i \le x_{i+1} + (1-y_i)      &(i\in T \setminus \{u\}) &\mbox{(C3)}\\
&&x_i \ge x_{i+1} - y_i          &(i\in T \setminus \{u\}) &\mbox{(C4)}\\
&&y_i \ge y_{i+1}                &(i \in T \setminus \{u\}) &\mbox{(C5)}\\
&&y_{i} \in \{0,1\}              &(i\in T)
\end{array}
\end{align*}
where $d(X,P)$ is some distance metric between $X$ and $P$. 
(C1) and (C2) ensure that $X$ is a distribution.
Each $y_i$ is a binary variable indicating whether index $i$ is before the peak. Monotonicity and single-peakedness are ensured by (C3) and (C4), and (C5), respectively. 
We propose to use the mean squared error (MSE)
\begin{align*}
\frac{1}{|T|}\sum_{i \in T} (x_i -p_i)^2
\end{align*}
for the objective function to ease computation. 
We also tested the Earth mover's distance and the mean absolute error (MAE); however, these did not demonstrate satisfactory performance. 
In the evaluation, we use the Jensen–Shannon divergence (JSD), a symmetrized version of the Kullback–Leibler divergence (KLD), defined as follows: 
\begin{align*}
D_{\rm KL}(P_1,P_2)&=\sum_{i \in T} p_{1i} \log 
\frac{p_{1i}}{p_{2i}}
,&\\
D_{\rm JS}(P_1,P_2)&=\frac{1}{2}D_{\rm KL}(P_1,(P_1+P_2)/2)+ \frac{1}{2}D_{\rm KL}(P_2,(P_1+P_2)/2).
\end{align*}

\subsection{Unimodal regression with a stochastic order constraint}
\label{sec:ours}
Suppose we estimate the search timing distributions $X_1$ and $X_2$ for queries \texttt{Q1} and \texttt{Q2}. Given prior knowledge that \texttt{Q1} should \emph{precede} \texttt{Q2}, we impose $X_1 \le_{\rm st} X_2$, yielding the following formulation:
\begin{align*}
\begin{array}{lllll}
R_{\rm st}(P_1,P_2):
&\mbox{min.}    &\displaystyle d(X_1,P_1)+d(X_2,P_2)\\\\
&\mbox{s. t.}  &\displaystyle \sum_{i\in T} x_{1i} = 1,
&\displaystyle \sum_{i\in T} x_{2i} = 1&\\
&&x_{1i} \in [0,1], &x_{2i} \in [0,1]                                       &(i\in T)\\
&&x_{1i} \le x_{1,i+1} + 1-y_{1i}, &x_{2i} \le x_{2,i+1} + 1-y_{2i}     &(i\in T \setminus \{u\})\\
&&x_{1i} \ge x_{1,i+1} - y_{1i}, &x_{2i} \ge x_{2,i+1} - y_{2i}             &(i\in T \setminus \{u\})\\
&&y_{1i} \ge y_{1,i+1},  &y_{2i} \ge y_{2,i+1}                              &(i\in T \setminus \{u\})\\
&&y_{1i} \in \{0,1\}, &y_{2i} \in \{0,1\}                &(i\in T)\\\\
&&\displaystyle \sum_{i \le t} x_{1i} \ge \sum_{i \le t} x_{2i}     &&(t \in T)
\end{array}
\end{align*}
where $P_1$ and $P_2$ are the empirical discrete distributions, and the last constraint imposes $X_1 \le_{\rm st} X_2$ as indicated in $(\ref{st_d})$. 
The following points should be noted: 
\begin{itemize}
\item In our application, we use $\le_{\rm st}$ because our prior knowledge concerns only the location of the distributions. However, our model readily extends to alternatives such as the convex and increasing convex orders, since their definitions also admit linear representations;
\item With sufficient data, the validity of imposing the stochastic order can be statistically tested prior to estimation; see \cite{lando2025new,barrett2003consistent} for the corresponding testing procedures. \end{itemize}

\section{Proof of concept}
\label{sec:poc}

\begin{figure}[tb]
\centering
\includegraphics[clip, width=0.7\textwidth]{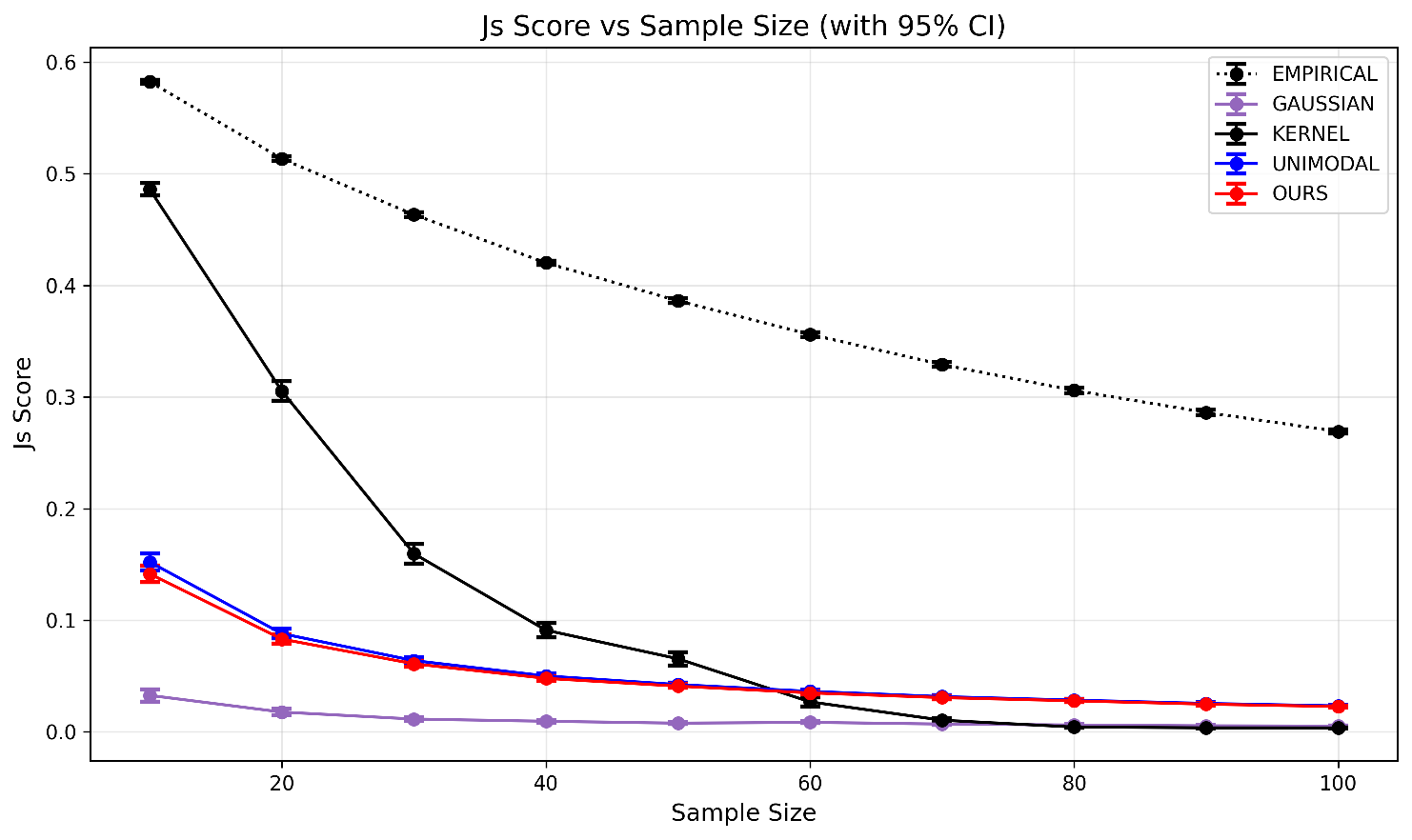}
\caption{Change in the estimation error (JSD) with respect to sample size using synthetic data}\label{fig:poc_sample}
\end{figure}

We first validate our model using a toy example. 
For comparison, we evaluate the following five methods:
\begin{itemize}
\item \texttt{EMPIRICAL}: the empirical distribution;
\item \texttt{GAUSSIAN}: the Gaussian maximum likelihood estimator;
\item \texttt{KERNEL}: a continuous distribution fitted via kernel density estimation using \texttt{KDEMultivariate} from \texttt{statsmodels} (ver.~0.13.1), with the bandwidth selected from ${0.05, 0.10, \ldots, 0.95}$ by two-fold cross-validation to minimize the mean squared error;
\item \texttt{UNIMODAL}: a unimodal distribution estimated by the model in Section~\ref{sec:unimodal};
\item \texttt{OURS}: multiple unimodal distributions estimated under stochastic order constraints using the model in Section~\ref{sec:ours}.
\end{itemize}
All code was implemented in Python, and the optimization problems were solved using 
\texttt{Gurobi Optimizer}\footnote{\url{https://www.gurobi.com/products/gurobi-optimizer/}} (ver. 9.0.3) on a machine running 
    macOS Sonoma 14.1.2 (Apple M1, 8GB). 

Suppose that the search timing distributions 
$X_1$ and $X_2$ for two search queries are given by normal distributions $N(-20,50)$ and $N(20,50)$. 
It is known~\cite{muller2001stochastic} that two normal distributions with identical variances satisfy the stochastic order if and only if their means are ordered; thus, $X_1 \le_{st} X_2$.
We constructed their empirical distributions, $P_1$ and $P_2$, by randomly generating $n$ samples for each query, where the sample size $n$ ranges from $10$ to $100$.

Figure~\ref{fig:poc_sample} illustrates the change in the estimation error (JSD) with respect to the sample size, where the bars indicate the confidence intervals computed from 100 trials. 
Since the data are generated from normal distributions, \texttt{GAUSSIAN} naturally achieves the lowest error. 
When the sample size is small ($n < 40$), \texttt{OURS} considerably outperformed \texttt{KERNEL}, but only slightly outperformed the conventional one, \texttt{UNIMODAL}. 
As the sample size increases, the accuracy of \texttt{KERNEL} improves, and the JSD scores of all methods eventually converge to a lower value.

\begin{figure}[tb]
\centering
\includegraphics[clip, width=\textwidth]{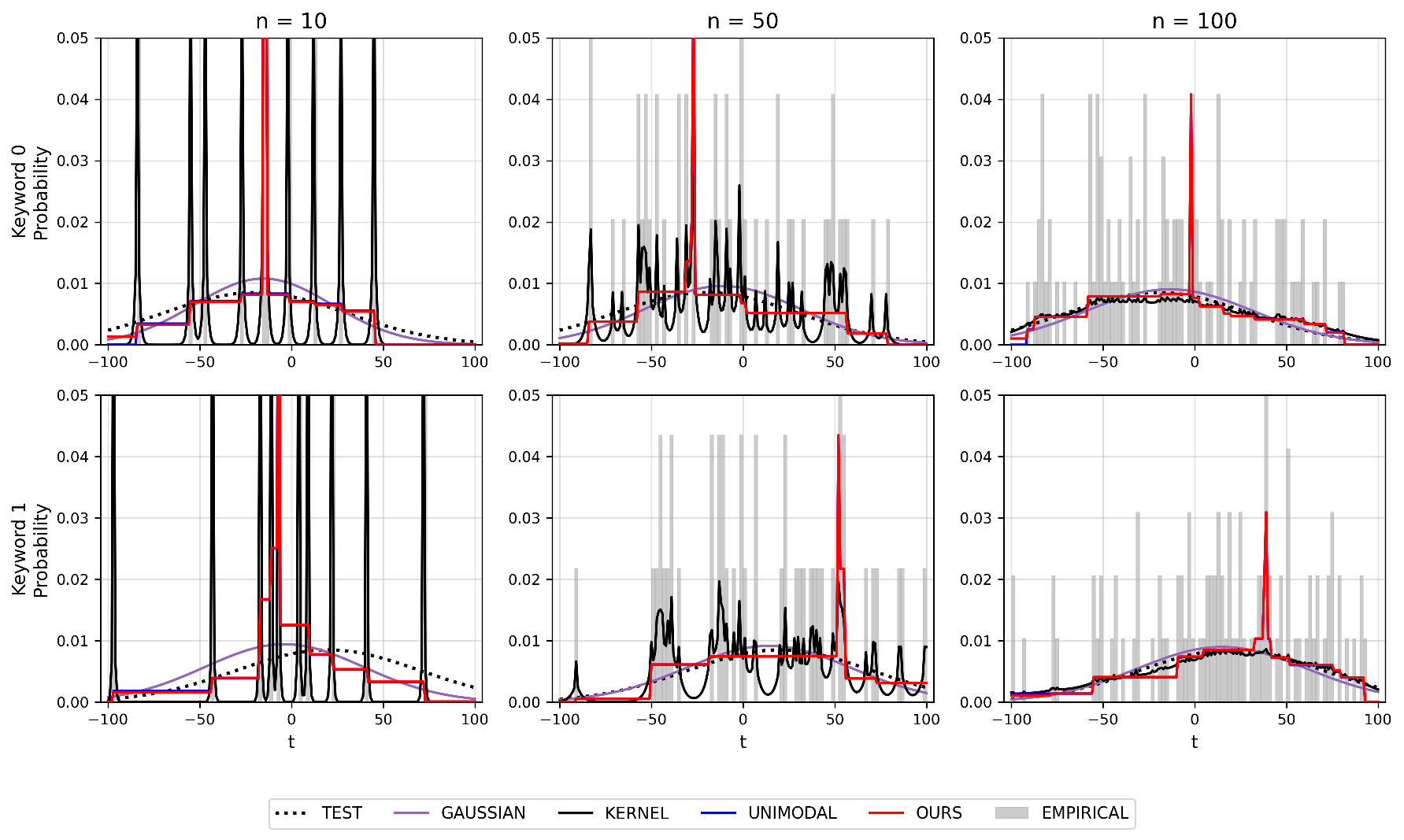}
\caption{Visual comparison of estimated distributions with varying sample sizes}\label{fig:poc_dist}
\end{figure}

To better understand these behaviors, Figure~\ref{fig:poc_dist} provides a visual comparison of the estimated distributions. 
In this figure, the upper row presents the estimation results for $X_1$, whose true search timing distribution follows a normal distribution $N(-20,50)$, while the lower row presents those for $X_2$, whose true distribution follows $N(20,50)$. Here, \texttt{TEST} corresponds to the true distribution. When the sample size is insufficient ($n=10$), \texttt{KERNEL} suffers from overfitting due to data sparsity, resulting in a spiky distribution that deviates from the true distribution.
In contrast, \texttt{OURS} successfully recovers the unimodal shape and the correct peak positions even with limited samples by leveraging the prior knowledge of unimodality and the stochastic order constraint.
This robustness against data sparsity leads to the superior JSD performance observed in Figure~\ref{fig:poc_sample}.

On the other hand, we see that \texttt{OURS} and \texttt{UNIMODAL} produce extremely steep distributions due to the monotonicity constraint. 
Such estimates may cause practical issues in applications, particularly when computing summary statistics such as quantiles. Therefore, in practice, it may be necessary to apply additional smoothing techniques to the estimated results or to incorporate appropriate regularization at the estimation stage, such as imposing constraints or penalties that discourage large differences between adjacent variables $x_{i}$ and $x_{i+1}$.

\section{Numerical Experiment}
\label{sec:ne}

Next, we report the experimental results on a real-world dataset accumulated on Mamari and provided by Connehito Inc. The computational environment and comparison methods are the same as those in the previous section.

\subsection{Setting}
This dataset originally consisted of 96,662,493 records collected from January 1, 2021, to December 31, 2022. 
We preprocessed this dataset as follows. 
\begin{enumerate}
\item For simplicity, we extracted users who have registered only one child, totaling 407,849 users. 
\item We then extracted the queries searched by at least one of the users identified in (1) during the specified period, totaling 29,836,147 queries.
\end{enumerate}
We define an \emph{instance} as a set of search queries in the form $(\texttt{Q1},\texttt{Q2},\ldots,\texttt{Qk})$, where each search query is a pair of keywords such as, \texttt{language\_1\_month\_old}. 
We considered the following three forms:  
\begin{enumerate}
\item[(F1)] $(
\texttt{X\_first\_trimester}, \,
\texttt{X\_second\_trimester}, \,
\texttt{X\_third\_trimester})$
\item[(F2)] $(
\texttt{X\_1\_month\_old}, \,
\texttt{X\_2\_months\_old}, \, \ldots, \,
\texttt{X\_6\_months\_old})$
\item[(F3)] $(
\texttt{X\_1\_year\_old}, \,
\texttt{X\_2\_years\_old}, \, \ldots, \,
\texttt{X\_12\_years\_old})$
\end{enumerate}
where \texttt{X} represents a specific term. 
Examples of \texttt{X} in (F1) include \texttt{diarrhea} and \texttt{body weight}. 
We selected certain terms as \texttt{X} if the corresponding search query yielded 80 or more records from January 2021 to December 2021 (training) and January 2022 to December 2022 (testing). 
As a result, we obtained 27 instances, as listed in the first column of Table~\ref{table:model_comparison_10}.
The number of records per search query ranged from 80 to 28,276.
We constructed datasets by randomly sampling a fixed number of records $n\in \{10,20,\ldots,80\}$ for each search query.

\subsection{Results}

\begin{table}[t]
\caption{Average estimation error (JSD) for each model with varying record sizes}
\label{tbl:avg_jsd_real}
\begin{tabular}{lrrrrr}
\toprule
Data set
& \texttt{EMP}
& \texttt{GAUSSIAN}
& \texttt{KERNEL}
& \texttt{UNIMODAL}
& \texttt{OURS}
\\
\midrule
10 records & 0.1425 & 0.0932 & 0.0983 & \underline{0.0888} & \textbf{0.0869} \\
20 records & 0.0859 & 0.1043 & \underline{0.0568} & 0.0571 & \textbf{0.0561} \\
30 records & 0.0656 & 0.1147 & \underline{0.0450} & 0.0453 & \textbf{0.0446} \\
40 records & 0.0552 & 0.1218 & \textbf{0.0375} & 0.0391 & \underline{0.0386} \\
50 records & 0.0487 & 0.1379 & 0.0354 & \underline{0.0350} & \textbf{0.0348} \\
60 records & 0.0442 & 0.1261 & \textbf{0.0310} & 0.0323 & \underline{0.0322} \\
70 records & 0.0410 & 0.1379 & \underline{0.0304} & \textbf{0.0302} & \textbf{0.0302} \\
80 records & 0.0385 & 0.1489 & 0.0302 & \textbf{0.0286} & \underline{0.0287} \\
\bottomrule
\end{tabular}
\end{table}

Table~\ref{tbl:avg_jsd_real} reports the estimation errors (JSD) of all methods, averaged over the 27 instances for each number of records. 
For each instance, the smallest estimation error is highlighted in bold, and the second smallest is indicated with an underline. 
As observed in the previous section, \texttt{OURS} considerably outperforms the baselines when the number of records is small ($n<40$), and as the number of records increases, this advantage diminishes. Unlike in the previous section, even when the number of records is large, \texttt{KERNEL} does not outperform \texttt{OURS}; rather, their performances are comparable, and in some cases \texttt{OURS} performs slightly better.
Since \texttt{GAUSSIAN} does not achieve good performance for this dataset, its results are omitted from the subsequent detailed experiments


\begin{table}[tb]
\caption{Estimation error by each model 
for the 10-record dataset}
\label{table:model_comparison_10}
\begin{tabular}{clrrrr}
\toprule
\multicolumn{2}{l}{Instance}&\texttt{EMP}&\texttt{KERNEL}&\texttt{UNIMODAL}&\texttt{OURS}\\
\midrule
(F1)&\texttt{diarrhea} & 0.238 & 0.157 & \underline{0.115} & \textbf{0.113} \\
&\texttt{body\_weight} & 0.251 & 0.171 & \underline{0.131} & \textbf{0.130} \\
&\texttt{constipation} & 0.209 & 0.128 & \underline{0.097} & \textbf{0.095} \\
&\texttt{abdominal\_pain} & 0.205 & 0.135 & \underline{0.111} & \textbf{0.109} \\
&\texttt{anemia} & 0.213 & 0.139 & \underline{0.109} & \textbf{0.108} \\\midrule
(F2)&\texttt{toddler\_defiance} & 0.122 & \underline{0.089} & 0.090 & \textbf{0.087} \\
&\texttt{not\_walking} & 0.195 & \underline{0.139} & \textbf{0.113} & \textbf{0.113} \\
&\texttt{language} & 0.130 & 0.092 & \underline{0.086} & \textbf{0.080} \\\midrule
(F3)&$\emptyset$
& 0.104 & \underline{0.074} & 0.075 & \textbf{0.071} \\
&\texttt{lying\_down\_on\_stomach} & 0.137 & 0.099 & \underline{0.091} & \textbf{0.089} \\
&\texttt{self-settle} & 0.126 & 0.086 & \textbf{0.082} & \underline{0.083} \\
&\texttt{body\_weight} & 0.128 & 0.090 & \underline{0.085} & \textbf{0.082} \\
&\texttt{not\_gaining\_weight} & 0.139 & 0.097 & \underline{0.092} & \textbf{0.090} \\
&\texttt{crying\_at\_night} & 0.108 & \underline{0.076} & \underline{0.076} & \textbf{0.074} \\
&\texttt{not\_sleeping} & 0.106 & \underline{0.077} & \underline{0.077} & \textbf{0.075}\\
&\texttt{not\_turning\_over} & 0.186 & 0.133 & \underline{0.109} & \textbf{0.104} \\
&\texttt{thumb\_sucking} & \underline{0.109} & \textbf{0.077} & \textbf{0.077} & \textbf{0.077} \\
&\texttt{breastfeeding\_interval} & 0.113 & 0.079 & \underline{0.078} & \textbf{0.075} \\
&\texttt{formula\_milk\_interval} & 0.114 & 0.081 & \underline{0.080} & \textbf{0.079} \\
&\texttt{exclusive\_breastfeeding\_interval} \quad & 0.112 & 0.079 & \underline{0.076} & \textbf{0.074} \\
&\texttt{mix\_feeding\_interval} & 0.109 & \underline{0.077} & \underline{0.077} & \textbf{0.075} \\
&\texttt{daily\_rhythm} & 0.126 & 0.089 & \underline{0.086} & \textbf{0.082} \\
&\texttt{sleep\_duration} & 0.114 & 0.079 & \underline{0.078} & \textbf{0.076}\\
&\texttt{sleep\_regression} & 0.112 & \underline{0.079} & \textbf{0.078} & \textbf{0.078} \\
&\texttt{comfort\_nursing} & 0.104 & \textbf{0.073} & \underline{0.074} & \underline{0.074} \\
&\texttt{weaning\_amount} & 0.119 & 0.081 & \underline{0.078} & \textbf{0.076} \\
&\texttt{not\_eating\_weaning} & 0.119 & \underline{0.080} & \underline{0.080} & \textbf{0.077} \\
\bottomrule
\end{tabular}
\end{table}

Table~\ref{table:model_comparison_10} presents the estimation error of each model for each instance in the 10-record dataset. 
Here, $\emptyset$ in the first row of (F3) represents an instance in which \texttt{X} is empty, corresponding to cases where users search for terms such as \texttt{1\_year\_old} alone, without combining them with other terms. 
We see from this table that \texttt{OURS} achieves the best performance for all but two instances; even in those two cases, it attains the second-best performance.  
Similarly, \texttt{UNIMODAL} achieves either the best or the second-best performance in all but two instances. 
Compared to \texttt{EMP}, the reduction rate of the estimation error by \texttt{OURS} was 
$36.87\%$ on average, $54.29\%$ at most, and $28.34\%$ at least.
Compared to \texttt{KERNEL}, the reduction rate was 
$9.31\%$ on average, $27.97\%$ at most, and $-1.65\%$ at least.
Finally, compared to \texttt{UNIMODAL}, the reduction rate was even smaller, 
$2.19\%$ on average, $6.35\%$ at most, and $-0.68\%$ at least.


\begin{table}[tb]
\centering
\caption{Estimation error by each model 
for the 80-record data set}
\begin{tabular}{clrrrr}
\toprule
\multicolumn{2}{l}{Instance}&\texttt{EMP}&\texttt{KERNEL}&\texttt{UNIMODAL}&\texttt{OURS}\\
\midrule
(F1)&\texttt{diarrhea} & 0.068 & \underline{0.041} & \textbf{0.037} & \underline{0.041} \\
&\texttt{body\_weight} & 0.075 & \underline{0.061} & \textbf{0.060} & \underline{0.061} \\
&\texttt{constipation} & 0.044 & \underline{0.033} & \textbf{0.030} & \textbf{0.030} \\
&\texttt{abdominal\_pain} & 0.058 & \textbf{0.037} & \textbf{0.037} & \underline{0.038} \\
&\texttt{anemia} & 0.051 & 0.037 & \textbf{0.036} & \textbf{0.036} \\\midrule
(F2)&\texttt{toddler\_defiance} & 0.039 & \textbf{0.030} & 0.032 & \underline{0.031} \\
&\texttt{not\_walking} & 0.068 & \underline{0.048} & \textbf{0.044} & 0.049 \\
&\texttt{language} & 0.039 & 0.031 & \underline{0.028} & \textbf{0.027} \\\midrule
(F3)&$\emptyset$ & 0.027 & 0.022 & \underline{0.019} & \textbf{0.017} \\
&\texttt{lying\_down\_on\_stomach} & 0.045 & \underline{0.040} & \textbf{0.036} &\textbf{0.036} \\
&\texttt{self-settle} & 0.030 & \textbf{0.023} & \underline{0.024} & 0.025 \\
&\texttt{body\_weight} & 0.037 & 0.032 & \underline{0.031} & \textbf{0.030} \\
&\texttt{not\_gaining\_weight} & 0.044 & \textbf{0.033} & \underline{0.034} & \underline{0.034} \\
&\texttt{crying\_at\_night} & 0.028 & \underline{0.024} & \textbf{0.022} & \textbf{0.022} \\
&\texttt{not\_sleeping} & 0.024 & 0.025 & \underline{0.021} & \textbf{0.020} \\
&\texttt{not\_turning\_over} & 0.069 & \underline{0.049} & \textbf{0.046} & \textbf{0.046} \\
&\texttt{thumb\_sucking} & 0.025 & \underline{0.023} & \textbf{0.020} & \textbf{0.020} \\
&\texttt{breastfeeding\_interval} & 0.026 & 0.022 & \underline{0.019} & \textbf{0.018} \\
&\texttt{formula\_milk\_interval} & 0.027 & \underline{0.024} & \textbf{0.023} & \textbf{0.023} \\
&\texttt{exclusive\_breastfeeding\_interval} \quad & 0.025 & \underline{0.024} & \textbf{0.020} & \textbf{0.020} \\
&\texttt{mix\_feeding\_interval} & 0.024 & \underline{0.023} & \textbf{0.020} & \textbf{0.020} \\
&\texttt{daily\_rhythm} & 0.034 & \underline{0.028} & \underline{0.028} & \textbf{0.027} \\
&\texttt{sleep\_duration} & 0.030 & 0.025 & \underline{0.024} & \textbf{0.023} \\
&\texttt{sleep\_regression} & 0.023 & \underline{0.020} & \textbf{0.019} & \textbf{0.019} \\
&\texttt{comfort\_nursing} & 0.024 & \underline{0.023} & \textbf{0.021} & \textbf{0.021} \\
&\texttt{weaning\_amount} & 0.027 & \underline{0.020} & \underline{0.020} & \textbf{0.019} \\
&\texttt{not\_eating\_weaning} & 0.030 & \textbf{0.020} & \underline{0.022} & \underline{0.022} \\
\bottomrule
\end{tabular}
\label{table:model_comparison_80}
\end{table}

Table~\ref{table:model_comparison_80} is the 
80-record dataset version
of Table~\ref{table:model_comparison_10}.
While \texttt{OURS} has a smaller advantage over the baselines than in Table~\ref{table:model_comparison_10}, it still outperforms them in most instances. 
Specifically, out of the 27 instances, \texttt{OURS} achieves the smallest estimation error in 19 instances and attains either the smallest or the second-smallest estimation error in 25 instances; similarly, \texttt{UNIMODAL} achieves the smallest estimation error in 15 instances and attains either the smallest or the second-smallest estimation error in 26 instances. 
Compared to \texttt{EMP}, \texttt{KERNEL}, and \texttt{UNIMODAL}, 
the reduction rates in the estimation error were 
$24.07\%$,  
$5.75\%$, and 
$0.35\%$ on average  
ranging from $12.37\%$, $-9.66\%$, and $-11.70\%$ at least to
$39.59\%$, 
$20.49\%$, and 
$7.68\%$ at most, respectively. 

\subsection{Discussion}

In Table~\ref{tbl:avg_jsd_real}, when the number of records is small, \texttt{OURS} outperforms \texttt{KERNEL} and \texttt{UNIMODAL}. This improvement can be attributed to the unimodality and stochastic order constraints, which allow the estimator to effectively pool the limited samples allocated to each distribution and estimate them jointly. Indeed, in Tables~\ref{table:model_comparison_10} and~\ref{table:model_comparison_80}, we confirm that the cases in which \texttt{OURS} achieves relatively better performance tend to involve a larger number of distributions estimated simultaneously. On the other hand, as the number of records increases, the unimodality and stochastic order constraints are often satisfied automatically. In such cases, imposing these constraints may unnecessarily restrict the flexibility of the estimator, which can in turn deteriorate estimation accuracy.

Compared with the previous section, Table~\ref{table:model_comparison_10} shows that \texttt{GAUSSIAN} does not perform well. One possible reason is that the target distributions are asymmetric and steep, and thus deviate substantially from normality. This may also explain why the performance of \texttt{KERNEL} does not improve significantly as the number of records increases, since kernel-based methods generally require the underlying distribution to be reasonably smooth in order to perform well.

\begin{figure}[tb]
\subfigure[\texttt{body\_weight} in (F1)]{%
\includegraphics[clip, width=0.7\textwidth]{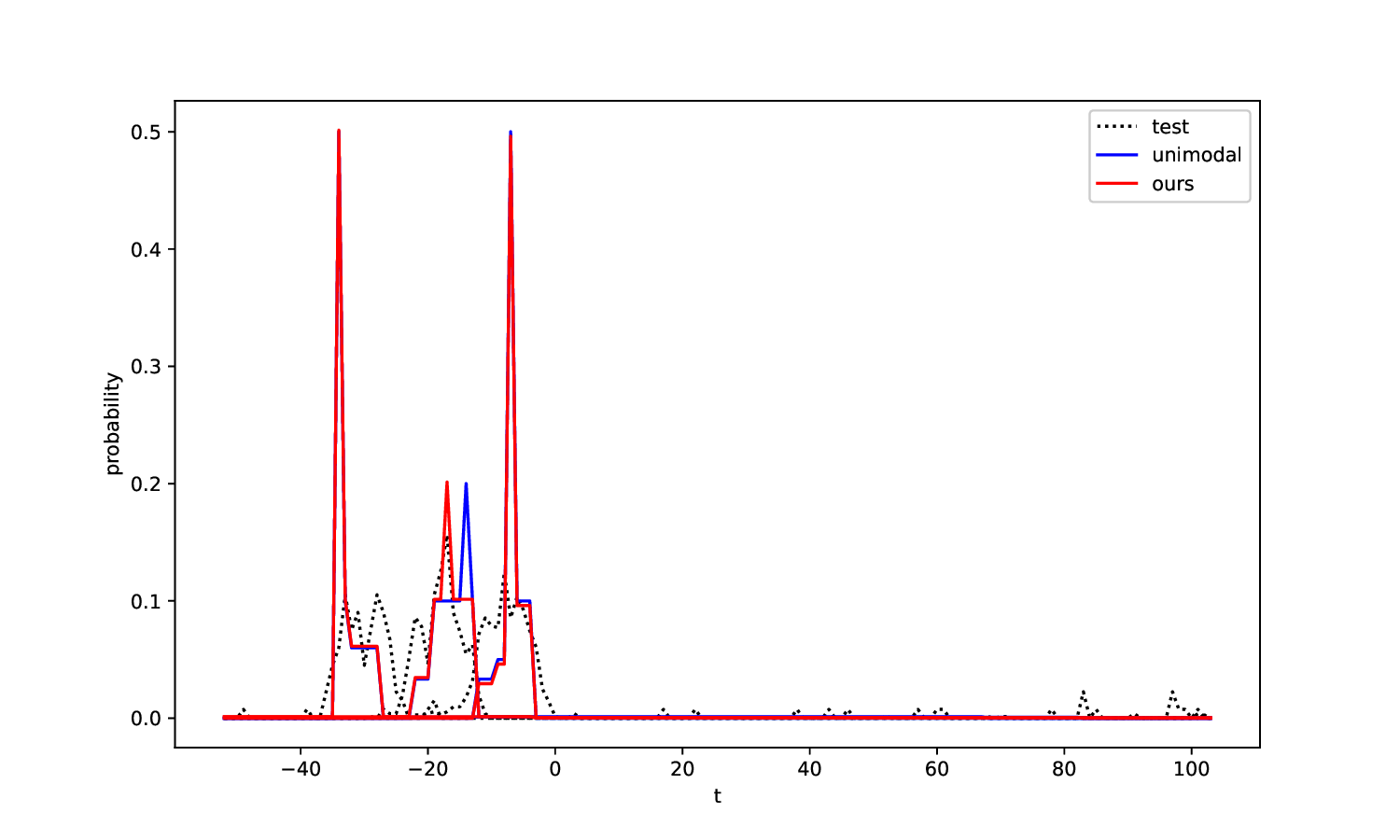}}\\%
\subfigure[\texttt{diarrhea} in (F1)]{%
\includegraphics[clip, width=0.7\textwidth]{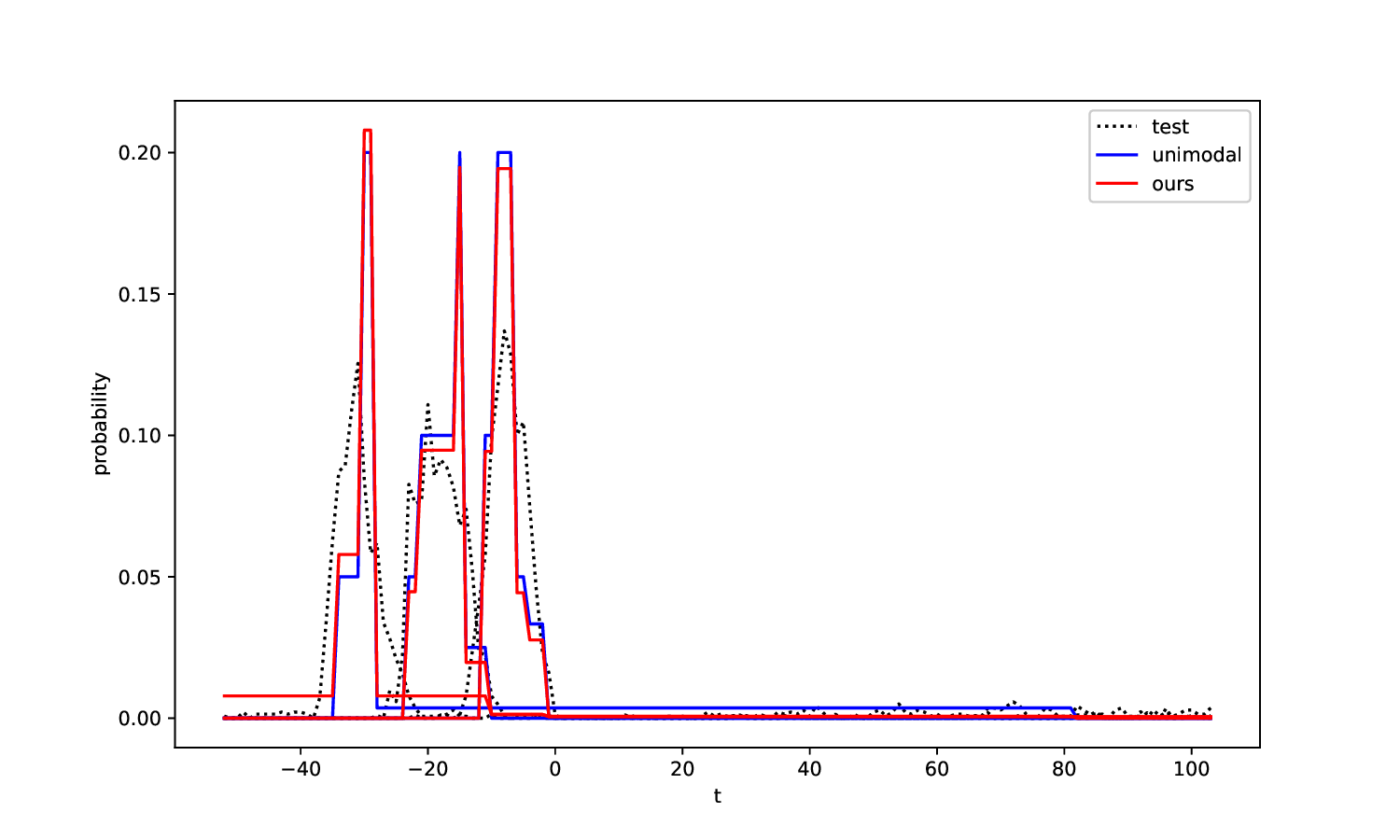}}%
\caption{Estimated distributions for two instances}
\label{fig:good_bad}
\end{figure}

We next examine cases where \texttt{UNIMODAL} and \texttt{OURS} yield different performance. Figure~\ref{fig:good_bad} presents two representative examples. In each figure, dotted lines indicate the true distributions used for evaluation, while solid lines show the estimated distributions.
In (a), the three distributions correspond to the first, second, and third trimesters. 
The left and right distributions are nearly identical for \texttt{UNIMODAL} (blue) and \texttt{OURS} (red). However, for the middle distribution (second trimester), the blue curve is slightly shifted to the left relative to the red curve due to the stochastic order constraint imposed in \texttt{OURS}. This adjustment reduces the estimation error of \texttt{OURS} compared with \texttt{UNIMODAL}. 
In (b), the red curve for the left distribution exhibits a heavier left tail induced by the stochastic order constraint, resulting in a noticeable deviation from the blue curve. In this case, the constraint adversely affects accuracy, and \texttt{OURS} yields a larger estimation error than \texttt{UNIMODAL}. 

Finally, in terms of computation time, \texttt{KERNEL} and \texttt{OURS} required a few seconds for all datasets, while \texttt{UNIMODAL} terminated in less than one second.

\section{Conclusion}
\label{sec:conc}

We addressed the problem of estimating the timing of specific keyword searches on an information site. We proposed a mixed-integer convex quadratic optimization model that incorporates prior knowledge of precedence relations among search timing distributions and can be efficiently solved by solvers. Through experiments on both synthetic and real-world datasets, we confirmed that our model improves estimation accuracy when the sample size is small. In terms of JSD, the proposed method achieves an average reduction of 2.2\%, with a maximum improvement of 6.3\% and a maximum deterioration of 0.7\%. When the sample size is sufficiently large, the proposed method performs comparably to existing models. 
Future work includes:
\begin{itemize}
\item investigating broader application settings involving multiple distributions with natural precedence relations, such as marketing analyses that track changes in customer interest;
\item developing methods to automatically determine which stochastic order constraints should be imposed, beyond the standard one considered in this study; and
\item 
advancing the theoretical aspects of the framework, including analysis of the properties of the estimators and the development of models and algorithms for producing smoother estimates.
\end{itemize}

\section*{Acknowledgment}
We would like to thank Connehito Inc. for providing the valuable dataset. We are also grateful to the three anonymous reviewers for their constructive and insightful comments, which led to substantial improvements in the numerical experiments.

\bibliography{mybib}
\bibliographystyle{abbrv}
\end{document}